\documentstyle[11pt,twoside]{article}

\input{mssymb}

\makeatother

\newcommand{\N}{{\Bbb N}}
\newcommand{\R}{{\Bbb R}}

\newcommand{\C}{{\Bbb C}}

\newcommand{\J}{{\cal J}}
\newcommand{\K}{{\cal K}}
\newcommand{\LL}{{\cal L}}

\newcommand{\Id}{I\mkern-1mud}
\newcommand{\eps}{\varepsilon}
\renewcommand{\epsilon}{\varepsilon}
\newcommand{\bea}{\begin{eqnarray*}}
\newcommand{\eea}{\end{eqnarray*}}
\newcommand{\beq}{\begin{equation}}
\newcommand{\eeq}{\end{equation}}
\newcommand{\rest}[2]{#1\raisebox{-0.5ex}{\mbox{$\mid_{#2}$}}} 
\newcommand{\Lp}{L_p}
\newcommand{\lp}{\ell_p}
\newcommand{\llq}{\ell_q}
\newcommand{\lplqn}{\ell_p^{\vphantom{n}}(\ell_q^n)}
\newcommand{\iy}{\infty}
\newcommand{\dopu}{{:}\allowbreak\ }

\newcommand{\pel}{Pe{\l}czy\'{n}ski}
\newcommand{\linq}{\mathop{\overline{{\rm lin}}}\nolimits}

\newcommand{\mysec}[2]{\begin{center}{\sc {#1}. {#2}}\end{center}%
\nopagebreak\setcounter{theo}{0}\setcounter{section}{#1}}

\newcommand{\Proof}{\par\noindent{\em Proof. }}
\newcommand{\eop}{\nopagebreak\hspace*{\fill}$\Box$}

\newtheorem{theo}{Theorem}[section]
\newtheorem{lemma}[theo]{Lemma}

\newtheorem{prop}[theo]{Proposition}

\newcounter{abc}   
\newcounter{iiiii} 

{\end{list}}

\newenvironment{statements}%
{\setcounter{abc}{0}
\begin{list}%
{{\rm (\alph{abc})}}
{\usecounter{abc}
\parsep=0pt plus 1pt
\topsep=1pt plus 2pt minus 1pt
\itemsep=1pt plus 2pt minus 1pt
\leftmargin=3\baselineskip
\labelsep=.6\baselineskip
\labelwidth=2.4\baselineskip
\rightmargin 0pt}%
}
{\end{list}}

{\begin{list}%
{$\bullet$}
{\leftmargin=5\baselineskip
\labelsep=\baselineskip
\labelwidth=2.5\baselineskip
\parsep=0pt plus 1pt
\topsep=1pt plus 2pt minus 1pt
\itemsep=1pt plus 2pt minus 1pt
\rightmargin 0pt}%
}
{\end{list}}

{\begin{list}%
{$\bullet$}
{\leftmargin 17mm
\labelsep 3mm
\labelwidth 12mm
\rightmargin 0pt}%
}
{\end{list}}

\newif\ifrefsc
\let\thebibliographyalt=\thebibliography                                %
\def\thebibliography#1                                                  %
 {\ifrefsc
 \begin{center}{\sc References}\end{center}\small
 \list                    %
 {[\arabic{enumi}]}{\settowidth\labelwidth{[#1]}\leftmargin\labelwidth  %
 \advance\leftmargin\labelsep                                           %
 \usecounter{enumi}}                                                    %
 \def\newblock{\hskip .11em plus .33em minus .07em}                     %
 \sloppy\clubpenalty4000\widowpenalty4000                               %
 \sfcode`\.=1000\relax                                                  %
 \else\thebibliographyalt{#1}\fi}                                         %

\refsctrue

\begin{document}

\begin{center}
{\Large\bf The $M$-ideal structure of some algebras of bounded linear
operators}
\\[33 pt]
{\sc Nigel J.~Kalton {\it and\/} Dirk Werner}
\end{center}
\bigskip
\begin{quote}
\small
{\sc Abstract.}
Let $1<p,\,q<\iy$.		
It is shown for complex scalars
that there are no nontrivial $M$-ideals in
$\LL(\Lp[0,1])$ if $p\neq 2$, 
and $\K(\lplqn)$ is the only nontrivial $M$-ideal in $\LL(\lplqn)$.
\end{quote}
\bigskip
%
%
%
%
\mysec{1}{Introduction}
A subspace $J$ of a Banach space $X$ is called an $M$-ideal if there
is an $\ell_1$-direct decomposition of the dual space $X^*$ of $X$
into the annihilator $J^\bot$ of $J$ and some subspace $V\subset X^*$:
$$
X^* = J^\bot \oplus_1 V;
$$
it is called nontrivial if $\{0\}\neq J\neq X$. This notion was 
introduced by Alfsen and Effros \cite{AlEf} and has proved useful in Banach
space geometry, approximation theory and harmonic analysis;
see \cite{HWW} for a detailed account.

A number of authors have studied the $M$-ideal structure in $\LL(X)$,
the space of bounded linear operators on a Banach space $X$,
with special emphasis on the question whether $\K(X)$, the subspace
of compact operators, is an $M$-ideal; see for instance
\cite{Beh6}, \cite{ChJo1}, \cite{ChJo2}, \cite{Fli},
\cite{GKS}, \cite{HaLi},
\cite{Kal-M}, \cite{Lim3}, \cite{OjaCR}, \cite{OjaDi}, \cite{PaWe},
\cite{SmWa2}, \cite{Dirk2}, \cite{WW3}
or Chapter~VI in \cite{HWW}. In particular we mention the facts that,
for a Hilbert space $H$, the $M$-ideals of $\LL(H)$ coincide with
its closed two-sided ideals \cite{SmWa1} 
and that, for a subspace $X$ of $\lp$,
$\K(X)$ is an $M$-ideal in $\LL(X)$ if and only if $X$ has the
metric compact approximation property \cite{ChJo1}.

In this paper we show that in many cases the ideal
of compact operators is the only candidate for an $M$-ideal
in $\LL(X)$. For $X=\lp$ this was done by Flinn \cite{Fli}.
For the function space $L_{p}=L_{p}[0,1]$ it has long been known
that the compact operators do not form an $M$-ideal if $p\neq2$
(\cite{Lim3} or, for another approach, \cite{OjaDi}); 
in \cite[p.~252]{HWW} the problem is posed to determine
the $M$-ideal structure of $\LL(\Lp)$ completely.
This is done in section~2 where we prove that there are no
nontrivial $M$-ideals in the algebra of operators on the complex space
$\Lp$, $p\neq 1,2,\iy$. In section~3 we study the $M$-ideals
in $\LL(\lplqn)$, for $1<p,\,q<\iy$. We show, for complex spaces again,
that here the compact operators form the only
$M$-ideal, thus proving a conjecture of Cho and Johnson \cite{ChJo2}.
(Here, as usual, 
$\lp(E_n)$ denotes the $\lp$-direct sum of the Banach spaces
$E_n$, i.e., $\lp(E_n)$ consists of all sequences of vectors
$x_n\in E_n$ such that $\|(x_n)\|=(\sum \|x_n\|^p)^{1/p} < \iy$.)

For our arguments we need the	notion of the (spatial)
numerical range of an operator
$T$ on a complex Banach space $E$. This is the set
$$
V(T)=\{x^*(Tx)\dopu \|x^*\|=\|x\|=x^*(x)=1\}.
$$
The operator $T$ is called hermitian if $V(T)\subset \R$. It is a 
well-known fact that an operator on $\Lp$ for $p\neq2$ is hermitian	
if and only if it is a multiplication operator, i.e., 
$Tf=hf$ for some real-valued $h\in L_\infty$. (See \cite{NMR1} and
\cite{NMR2} for details.)
\medskip

\noindent
{\bf Acknowledgement.} This work was  done while the second-named
author  was  visiting the University of Missouri,
Columbia. It is his pleasure to express his gratitude to all
those who made this stay possible.
\bigskip
%
%
%
%
\mysec{2}{$M$-ideals of operators on $\Lp[0,1]$}
To prove the main result of this section, Theorem~\ref{A2},
we will  need the following lemmas.
%
%
%
%
\begin{lemma}\label{A1}
For all $\eps>0$ there is some $\delta>0$ such that whenever $A$ 
and $B$ are operators on a complex Banach space satisfying
$\|\sin A\| \le 1-\delta$, $\|\sin B\|\le 1$,
$\arcsin(\sin A)= A$, $\arcsin(\sin B)= B$ and
$\|\sin A - \sin B\| \le \delta$,
then $\|A-B\|\le\eps$.
The condition $\arcsin(\sin A)=A$ holds if $V(A)
\subset R := [-0.1, 1.1] \times [-0.1,0.1]\,i$.
\end{lemma}
\Proof
Consider the power series expansion
$$
\arcsin z = \sum_{k=1}^\infty c_k z^k, \qquad |z|\le 1,
$$
in which all $c_k\ge 0$. 
Since $\sin A$ is contractive, the series $\sum c_{k} (\sin A)^{k}$
converges, and the operator $\arcsin(\sin A)$ is defined (likewise
for $B$).

Write $v=\sin A$, $u=\sin B$, $w=u-v$. By the uniform continuity of
the $\arcsin$-function on $[0,1]$ we conclude for sufficiently
small $\delta$
\bea
\|B-A\| &=&
\| \arcsin (v+w) - \arcsin v \| \\
&=&
\left\| \sum_{k=1}^\infty c_k \left( (v+w)^k -v^k \right) \right\| \\
&\le&
\sum_{k=1}^\iy c_k \sum_{l=0}^{k-1} {k\choose l} \|v\|^l \|w\|^{k-l} \\
&=&
\sum_{k=1}^\iy c_k \left( (\|v\|+\|w\|)^k - \|v\|^k \right) \\
&=&
\arcsin(\|v\|+\|w\|) - \arcsin \|v\| \\
&\le& \eps
\eea
if $\|w\|=\|\sin A - \sin B\|\le \delta$.

We finally discuss the numerical range condition. If $V(A)\subset R$,
then  the spectrum $\sigma(A)$ is also contained in $R$ 
\cite[p.~19]{NMR1}. Now $\sin R$ is a subset of the open unit disk.
Hence
$\arcsin(\sin z) = z$ for $z\in \sigma(A)\subset R$, and we deduce
by the functional calculus that $\arcsin(\sin A) = A$.
\eop
%
%
%
%
\begin{lemma}\label{A1a}
Let $X$ be a Banach space, $\J \subset \LL(X)$  a two-sided
ideal and $P$ a projection  onto a complemented subspace
$E$ of $X$ which is isomorphic to $X$.
\begin{statements}
\item
If $P\in \J$, then $\J = \LL(X)$.
\item
If $E$ is $C$-isomorphic with $X$ and $\J$ contains an operator $T$
with $\|T-P\|< (C\|P\|)^{-1}$, then $\J=\LL(X)$.
\end{statements}
\end{lemma}
\Proof
(a) Let $j\dopu E \to X$ denote the canonical embedding and let 
$\Phi\dopu E \to X$ denote an isomorphism. Then
$\Id =( \Phi  P) P( j	\Phi^{-1}) \in \J$.

(b) We retain the above notation and let $S=(\Phi P)T(j\Phi^{-1})\in
\J$. Then, if $\|\Phi\|\,\|\Phi^{-1}\|\le C$,
\bea
\|S-\Id\| &\le& \|\Phi\| \, \|P\| \, \|T-P\| \, \|\Phi^{-1}\| \\
&\le& C \, \|P\| \, \|T-P\| ~<~ 1;
\eea
thus $S$ is invertible and $\Id= SS^{-1} \in \J$.
\eop
%
%
%
%
\begin{lemma}\label{A1b}
Let $1\le p<\iy$ and let $\varphi$ and $\varphi'$ be bimeasurable
bijective transformations on $[0,1]$ and $\lambda$ and $\lambda'$
be measurable functions on $[0,1]$. Suppose
$I\dopu f \mapsto \lambda \cdot f \circ \varphi$			and
$J\dopu f \mapsto \lambda' \cdot f \circ \varphi'$ 
are two isometric isomorphisms on $\Lp$ such that $\{\varphi \neq
\varphi'\}$ has positive measure. Then $\|I-J\|\ge 2^{1/p}$.
\end{lemma}
\Proof
It is enough to prove this for $J=\Id$, that is, for $\lambda'=1$ and
$\varphi'(\omega)=\omega$. By Lusin's theorem, $\{\omega\dopu 
\varphi(\omega)\neq\omega\}$ contains a compact set $C$ of positive
measure such that $\rest{\varphi}{C}$ is continuous. A compactness
argument now yields some compact subset $A\subset C$ of positive measure
such that $\varphi(A) \cap A = \emptyset$. Pick $\alpha$ so that
$f=\alpha \chi_{\varphi(A)}$ has norm~1. Then $\|If\|=1$, $If$ and $f$
have disjoint supports and thus $\|If-f\| = 2^{1/p}$.
\eop
\bigskip

We now state and prove the main theorem of this section.
%
%
%
%
\begin{theo}\label{A2}
If $1<p<\iy$ and $p\neq2$, then there are no nontrivial $M$-ideals
in $\LL(\Lp([0,1],\C))$.
\end{theo}
\Proof
Let $\{0\}\neq \J	\subset \LL(\Lp)$ be an $M$-ideal. Since $\Lp$ and
its dual are uniformly convex, $\J$ is a two-sided ideal \cite{ChJo2}. 
We
shall show that $\Id\in \J$, thus proving our claim that $\J=\LL(\Lp)$.

In fact, by Lemma~\ref{A1a} it is enough to show that some 
characteristic projection
$P_A \dopu f\mapsto \chi_A f$, for some Borel set of positive measure,
is in $\J$. This suffices because then
$\Lp(A)$ is isometrically isomorphic to $\Lp[0,1]$ \cite[p.~321]{Royden}. 

We now let $0<\eta <0.1$ be a small number. (It will become clear
in due course how small $\eta$ should actually be.) 
By \cite[Th.~V.5.4]{HWW}
$\J$ contains an operator $T$ with $\|T\|=1$ and
$V(T)\subset [-\eta, 1+\eta] \times [-\eta, \eta]\,i$. We are going 
to show that this `almost' hermitian operator $T$ is close to an
hermitian operator $f\mapsto hf$; and we will eventually prove
that $P_A\in \J$ for $A={\{h\ge \frac{1}{2}\}}$.

We first apply Theorem~4 from \cite[p.~28]{NMR1} to obtain that
$$
\left\|e^{itT}\right\| \le 1+\alpha |t| \qquad \forall |t|\le1,
$$
where $\alpha=\alpha(\eta) \to 0$ 
as $\eta \to 0$. That is, for small enough $\eta$
and $|t|\le 1$, the operators $e^{itT}$ are small bound
isomorphisms:
$$
\frac{1}{1+\alpha} \|f\|
\le \left\|e^{itT}f\right\| \le (1+\alpha)\|f\|.
$$
By a result due to Alspach \cite{Als2}, $e^{itT}$ is close to an isometric
isomorphism $I_t$:
$$
\left\|e^{itT}-I_t\right\| \le \beta,
$$
where $\beta=\beta(\eta) \to 0$ 
as $\alpha(\eta) \to 0$. But the isometric isomorphisms
on $\Lp$, $p\neq2$, are known to have the form
$$
I_t f = \lambda_t \cdot f \circ \varphi_t
$$
for some measurable function $\lambda_t$ and some 
bimeasurable bijection $\varphi_t$ \cite[p.~333]{Royden}.
Note that by Lemma~\ref{A1b}
$\|I_t - I_s\|\ge 2^{1/p}$ if $\{ \varphi_s \neq \varphi_t\}$
has positive measure.
This and the fact that $(e^{itT})_{t\in\R}$		is a uniformly continuous group
of operators show that  $\varphi_t = \varphi_0$ for $|t|\le1$
provided $\eta $ is small enough, and we must have that $\varphi_0
(\omega)=\omega$ since $e^{itT}=\Id$ for $t=0$. 
Furthermore, this enforces $|\lambda_t| = 1$ a.e. Writing
$\lambda_t = e^{ih_t}$, $-\pi\le h_t \le \pi$, and $M_t$	for the
multiplication operator $f\mapsto h_t f$ we finally obtain
$$
\left\|e^{itT} - e^{iM_t}\right\| \le \beta   
$$
for $|t|\le 1$.

Fix $t=\frac12$. By continuity of inversion we may
assume in addition that 
$$
\left\|e^{-iT/2} - e^{-iM_{1/2}}\right\| \le \beta,
$$
where $\beta =\beta(\eta)\to 0$ as $\eta \to 0$. Consequently
$$
\|\sin T/2 - \sin M_{1/2} \|\le \beta.
$$
If $B$ denotes the operator of multiplication with $\arcsin(\sin
h_{1/2})$, then $\sin B\allowbreak=\sin M_{1/2}$, $\arcsin(\sin B)=B$ and
$\|{\sin B}\|\le1$. Also
$$
\|\sin T/2\| \le \sum_{k=0}^{\iy}\frac{1}{(2k+1)!}
\left\|\frac T2\right\|^{2k+1}  \le {\textstyle \sinh \frac12} <1.
$$
According to Lemma~\ref{A1} we have for the multiplication operator
$M:=2B$, say $f \mapsto hf$,
$$
\|T-M\|\le \gamma,
$$
where $\gamma=\gamma(\eta)\to 0 $ 
as $\eta \to 0$. Since $\|T\|=1$, it follows that
$A=\{h\ge \frac{1}{2}\}$ has positive measure if $\eta$ is sufficiently
small.

Now consider $P_A \J P_A$, which can isometrically be identified with
a subspace of $\LL(\Lp(A))$. It is straightforward to check that it
is an $M$-ideal and thus a two-sided ideal. Moreover, since
$\|P_A T P_A - P_A M P_A\|\le \gamma $ and $\rest{M}{\Lp(A)}$ is
invertible with an inverse of norm $\le 2$, it follows that
$P_A T P_A$ is invertible in $\Lp(A)$ once $\gamma < \frac{1}{2}$.
Therefore the ideal $P_A\J P_A$ coincides with $\Lp(A)$.
This implies $P_A\in\J$, as requested.
\eop
\bigskip

In the case $p=2$ the $M$-ideals and the closed two-sided ideals
of $\LL(L_2)$ coincide, as already mentioned.
Therefore the ideal of compact operators is
the only nontrivial $M$-ideal in $\LL(L_2[0,1])$.
For the $M$-ideals in $\LL(L_1)$ and $\LL(L_\iy)$ see \cite{WW3}.

For future reference we record the following result that was 
established in the preceding proof.
%
%
%
%
\begin{prop}\label{A3}
Let $1\le p< \iy$, $p\neq 2$.
For all $\gamma>0$ there is some $\eta>0$ such that, whenever 
$X$ is a complex $\Lp(\mu)$-space and $T\in \LL(X)$ satisfies
$V(T) \subset [-\eta,1+\eta] \times [-\eta,\eta]\,i$, then there is 
some hermitian operator $M$ satisfying $\|T-M\|\le \gamma$.
\end{prop}
\bigskip
%
%
%
%
\mysec{3}{$M$-ideals of operators on $\lplqn$}
In this section we shall primarily deal with the space $\lplqn$,
${1<p,\,q<\iy}$.
It is known 
(e.g., \cite{Lim3}) that 
the compact operators form an $M$-ideal in $\LL(\lplqn)$, and in
\cite{ChJo2} Cho and Johnson exhibit an ideal in $\LL(\lplqn)$ that fails
to be an $M$-ideal. Actually, they conjecture that the compact operators
are the only $M$-ideal. Here we offer a proof of this conjecture
in the case of complex scalars.

The following lemma will turn out to be useful; for related results 
see \cite{CasKotL}.
%
%
%
%
\begin{lemma}\label{C1}
There is a constant $C$ such that, whenever $(k_n)$ is a sequence
of positive integers with $\limsup k_n =\iy$, then $\lplqn$ is 
$C$-isomorphic to $\lp^{{\vphantom{n}}}(\llq^{k_n})$.
\end{lemma}
\Proof
Let $E_1, E_2,\ldots$ be an enumeration of the $\llq^n$-spaces in which
each $\llq^n$-space is repeated infinitely often, and let
$X=\lp(E_n)$, $Y=\lp^{{\vphantom{n}}}(\llq^{k_n})$. 
It is enough to show that $X$ is
$C'$-isomorphic to $Y$ for some universal constant $C'$. This follows
immediately from an application of \pel's decomposition method
\cite[p.~54]{LiTz1}, since there are subspaces $U\subset X$, $V
\subset Y$ and natural isometric isomorphisms $X\cong \lp(X)$,
$X\cong Y \oplus_p U$, $Y\cong X \oplus V$. Note that $X$ and $V$
are in fact complemented by contractive projections, so $Y$ is 
$4$-isomorphic to $X\oplus_p V$. Thus the decomposition scheme
yields the desired isomorphism with $C'=16$.
\eop
\bigskip

And now for a technical lemma.
%
%
%
%
\begin{lemma}\label{C1a}
Let $1<p<\iy$ and let $(E_n)$ be a sequence of finite-dimensional
Banach spaces. Let $P_k$ denote the canonical finite-rank
projection from $\lp(E_n)$ onto $E_k$. Suppose $T\in\LL(\lp(E_n))$
satisfies $\sup_k \|P_k T P_k\|\le \eps$ for some $\eps>0$.
Then there is a subsequence $(k_n)$ of the positive integers such that
the inequality
$$\|PTP\|\le 3\eps $$
holds for the canonical projection $P$ from $\lp(E_n)$ onto $\lp(E_{k_n})$.
\end{lemma}
\Proof
We first establish two claims.
\begin{quote}
{\em Claim 1.}  For all $\beta>0$ and all $k\in \N$ there exists 
$m_1\in\N$ such that for all $m\ge m_1$
$$
\left\| \sum_{l\ge m} P_l T P_k \right\| \le \beta.
$$
\end{quote}
This holds since $(\sum_{l\ge m} P_l)_m$ is uniformly bounded,
converges to~0 pointwise and thus converges to~0 uniformly on the 
compact set $TP_k(B)$, where $B$ denotes the unit ball of $\lp(E_n)$.
\begin{quote}
{\em Claim 2.}  For all $\beta>0$ and all $k\in \N$ there exists 
$m_2\in\N$ such that for all $m\ge m_2$
$$
\|  P_k T P_m \| \le \beta.
$$
\end{quote}
This holds since $P_m^* \to 0$ pointwise and thus $\|P_k T P_m\| =
\|P_m^* T^* P_k^* \| \to 0$.

We now construct the desired sequence inductively. We put $k_1=1$.
Suppose $k_1,\ldots,k_n$ have already been constructed. We apply
Claim~1 with $k=k_n$, $\beta= \eps/2^n$, and we pick $k_{n+1}\ge m_1$
so that
$$
\sum_{j=1}^n \|P_{k_j} T P_{k_{n+1}}\|\le \frac{\eps}{2^{n+1}},
$$
which is possible by Claim~2. Then the resulting sequence $(k_n)$
satisfies
$$
\left\|\sum_{j>n} P_{k_j} T P_{k_{n}}\right\|\le
\left\|\sum_{l\ge k_{n+1}} P_{l} T P_{k_{n}}\right\|\le
\frac{\eps}{2^n}.
$$
Furthermore by assumption on $T$
$$
\left\|\sum_{n=1}^\iy P_{k_n} T P_{k_{n}}\right\| =
\sup_n \|P_{k_n} T P_{k_{n}}\| \le \eps,
$$
since $\sum_n P_{k_n} T P_{k_{n}}$ is a block-diagonal operator.
Therefore,
\bea
\|PTP\| &\le&
\sum_{n=1}^\iy \sum_{j<n} \|P_{k_j} T P_{k_{n}}\| +
\left\| \sum_{n=1}^\iy P_{k_n} T P_{k_{n}}\right\|	+
\sum_{n=1}^\iy \left\|\sum_{j>n} P_{k_j} T P_{k_{n}}\right\|	\\
&\le&
\sum_{n=1}^\iy \frac{\eps}{2^n} + \eps +
\sum_{n=1}^\iy \frac{\eps}{2^n}    \\
&=&3\eps.
\eea
\nopagebreak\eop
%
%
%
%
\begin{theo}\label{C2}
Consider the complex Banach space $\lplqn$, $1<p,\,q<\iy$. Then
$\K(\lplqn)$ is the only $M$-ideal in $\LL(\lplqn)$.
\end{theo}
\Proof
We have already pointed out that $\K(\lplqn)$ as a matter of fact is
an $M$-ideal. Moreover, it follows e.g.\ from \cite[Th.~4.4]{HaLi}
that an $M$-ideal $\J\neq\{0\}$ contains $\K(\lplqn)$, and since
$\lplqn$ and its dual are uniformly convex, we know from \cite{ChJo2} 
that an $M$-ideal in $\LL(\lplqn)$ must be a two-sided ideal.

Let us now assume that $\J$ is an $M$-ideal in $\LL(\lplqn)$ strictly
containing $\K(\lplqn)$. We have to prove that $\J=\LL(\lplqn)$.
We first claim that an operator $T\in \LL(\lplqn)$ factoring
through $\lp$ belongs to $\J$. In fact, let us write  $T=T_2 T_1$ with
$T_1\in \LL(\lplqn,\lp)$ and $T_2\in \LL(\lp,\lplqn)$. Now $\J$
contains a noncompact operator $S$. The proof of \cite[Prop.~2.c.3]{LiTz1}
shows that $S$ acts as an isomorphism on a complemented copy $E$ of 
$\lp$, and $F=S(E)$ is complemented, too. Let 
$\pi\dopu \lplqn\to E$ and $\sigma\dopu \lplqn\to F$ denote projections.
Then $\pi=(S^{-1}\sigma)S\pi\in \J$, and letting $\Phi\dopu\lp\to E$
denote an isomorphism, we see that $T=T_2T_1= T_2\Phi^{-1}\pi\Phi T_1
\in \J$.

Next we observe that it is enough to prove the theorem under the assumption
that $q\le p$; the remaining case $q>p$ then follows by duality. We may
also assume that $q\neq2$. Indeed, $\lp^{{\vphantom{n}}}(\ell_2^n)$ 
is known to be
isomorphic to a complemented subspace of $\lp$ and thus isomorphic to
$\lp$ itself. But the compact operators are the only closed 
two-sided ideal in $\LL(\lp)$ (see the above claim); hence
$\K(\lp^{{\vphantom{n}}}(\ell_2^n))$	
is the only $M$-ideal in $\LL(\lp^{{\vphantom{n}}}(\ell_2^n))$,
since every $M$-ideal is an ideal. So we will suppose that $q\le p$ and
$q\neq 2$ in the remainder of the proof.

We need some notation. We will denote the unit vectors in $\lplqn$
by $e_{kl}$, $k\in \N$, $1\le l\le n$. It is clear that 
$E_1= \linq \{e_{k1}\dopu k\in\N\}$ is isometric to $\lp$. If $q\le p$,
then the identity map from $\llq^n$ to $\lp^n$ is contractive; thus
there is a natural operator $A$ mapping $\lplqn$ to 
$\lp^{{\vphantom{n}}}(\lp^n)\cong \lp
\cong E_1$. We can think of this operator $A$ as an element of
$\LL(\lplqn)$, and by construction $A$ factors through $\lp$.	Consequently
$A\in\J$, $\|A\|=1$, and also $\|Ae_{kl}\|=1 $ for all $k$ and $l$.

At this stage we invoke a result from \cite{Dirk6} saying that there is a 
net $(H_\alpha)\subset \J$ converging to $\Id$ in
the $\sigma(\LL,\J^*)$-topology such that
\beq\label{eqC1}
\limsup \|{\pm A} + (\Id-H_\alpha)\| =1 .
\eeq
An application of \cite[Th.~V.5.4]{HWW} 
and a convex combinations technique
described there allows us
to assume in addition that the numerical ranges $V(H_\alpha)$ are 
contained in small rectangles 
$R_\alpha = [-\eta_\alpha,1+\eta_\alpha] 
\times [-\eta_\alpha,\eta_\alpha]\,i$,
with $\eta_\alpha\to 0$. We denote the canonical injection of $\llq^k$
 into $\lplqn$ by $j_k$ and the projection from $\lplqn$ onto $\llq^k$
by $P_k$. For each $k$, $V(P_k(\Id-H_\alpha)j_k)\subset R_\alpha$;
and therefore, by Proposition~\ref{A3}, 
we conclude 
that, given $\eps>0$, for large enough $\alpha$
all the operators $P_k(\Id-H_\alpha)j_k$, $k=1,2,
\ldots$, differ from hermitian operators $M_{\alpha k}\dopu \llq^k \to
\llq^k$ by less than $\eps/3$. Also, from (\ref{eqC1}) and the 
uniform convexity of $\lplqn$ we infer that $(\Id-H_\alpha)e_{kl}\to 0$
(since $\|Ae_{kl}\|=1$), uniformly in $k$ and $l$. From this we
deduce for large enough $\alpha$ that $\|M_{\alpha k}(e_{kl})\|\le
\frac{2}{3}\eps$ for all $k$ and $l$ and thus, since the $M_{\alpha k}$
are hermitian (i.e., multiplication operators), $\|M_{\alpha k}\|\le
\frac{2}{3}\eps$ for all $k$. 
This means for large enough $\alpha$ and $H=H_\alpha$
$$
\|P_k(\Id-H)P_k\|	\le \|P_k(\Id-H)j_k\| \le \eps
\qquad \forall k\in\N.
$$

Hence $\Id-H$ meets the assumption of Lemma~\ref{C1a}. That lemma
provides us with a sequence $(k_n)$ such that
for the canonical
projection $P$ from $\lplqn$ onto $\lp^{{\vphantom{n}}}(\llq^{k_n})$
$$
\|P-PHP\| = \|P(\Id-H)P\| \le 3\eps.
$$
Since $PHP\in\J$, an appeal to Lemmas \ref{A1a}(b) and \ref{C1} finishes
the proof of Theorem~\ref{C2}
provided we have chosen $\eps< \frac{1}{3C}$, where $C$ is 
the constant appearing in Lemma~\ref{C1}.
\eop
\bigskip

A similar, but technically simpler argument yields the following
result.
%
%
%
%
\begin{prop}\label{C4}
Let $1<p<\iy$, and let $X$ be a uniformly convex Banach space
isomorphic to a subspace of $\lp$. Suppose $\J\subset \LL(X)$
is a closed ideal strictly containing $\K(X)$. If $\J$ is an
$M$-ideal, then $\J=\LL(X)$.
\end{prop}
\Proof
As in the proof of Theorem~\ref{C2}, we argue that $\J$ contains all
$\lp$-factorable operators. Since $X$ embeds into $\lp$, $X$ contains a
subspace isomorphic to $\lp$. Let $\Psi\dopu X\to\lp$ and $\Phi\dopu
\lp\to X$ denote (into-)isomorphisms, and let $A=\Phi\Psi$.  Then
$A$ factors through $\lp$, and hence $A\in\J$.

We may assume that $\|A\|=1$. There is some $\eps>0$ such that 
$\eps\|x\| \le \|Ax\|$ for all $x\in X$. Again, there is a net
$(H_\alpha )\subset \J$ such that
$$
\limsup \|{\pm A} +(\Id-H_\alpha)\| =1.
$$
If $\alpha$ is large enough and $H=H_\alpha$, we have by uniform convexity
for some $\delta>0$ that $\|(\Id-H)x\| \le 1-\delta$ whenever $\|x\|=1$.
Hence $\|\Id-H\|<1$, and $\J$ contains an invertible element.
This proves that $\J=\LL(X)$.
\eop
\bigskip

We recall that $\K(X)$ is an $M$-ideal in $\LL(X)$ if $X$ is isometric
to a subspace of $\lp$ and has the metric compact approximation
property.
\bigskip

We finally mention some problems suggested by our work.
\smallskip

(1) The results in sections 2 and 3 support the conjecture that $\K(X)$
is the only candidate for an $M$-ideal in $\LL(X)$ if $X$ is 
isometric to a subspace of $\Lp$. Does this conjecture hold?

(2) In the other direction, it would be interesting to find an example
of a uniformly convex Banach space $X$ and an ideal $\K(X)\neq \J\subset
\LL(X)$ which is a nontrivial $M$-ideal. It is also still open whether,
for a uniformly convex Banach space $X$, every $M$-ideal in $\LL(X)$
is a two-sided ideal. (An $M$-ideal is known to be a left ideal
if $X$ is uniformly convex \cite{ChJo2}.)
\bigskip
%
%
%
%

%
%
%
%

\normalsize
\bigskip\bigskip\bigskip
\begin{tabbing}
\hspace*{.55 \textwidth} \= \kill
Department of Mathematics\> I.~Mathematisches Institut\\
University of Missouri  \> Freie Universit\"at Berlin \\
  \> Arnimallee 2--6 \\
Columbia, Mo 65\,211 \> 14\,195 Berlin \\[2pt]
U.S.A. \> Federal Republic of Germany			\\[4pt]
e-mail:\> e-mail: \\
mathnjk@mizzou1.missouri.edu  \>
werner@math.fu-berlin.de
\end{tabbing}

\end{document}